\documentclass{article}[12]

\usepackage{amsmath}
\usepackage{amsfonts}
\usepackage{amsthm}
\usepackage{graphicx}
\usepackage{float}

\theoremstyle{plain}
\newtheorem{thm}{Theorem}[section]
\newtheorem{lem}[thm]{Lemma}
\newtheorem{propo}[thm]{Proposition}

\theoremstyle{definition}
\newtheorem{defin}{Definition}[section]

\theoremstyle{remark}

\newcommand{\eqna}[1]{
    \begin{eqnarray}
        #1
    \end{eqnarray}
}
\newcommand{\eq}[1]{
    \begin{equation}
        #1
    \end{equation}
}

\newcommand{\Mm}{\mathcal{M}}
\newcommand{\hM}{\widehat{\mathcal{M}}_{g,1}}
\newcommand{\hD}{\widehat{D}}

\newcommand{\As}{\mathcal{A}}
\newcommand{\Fl}{\mathcal{F}}
\newcommand{\Mt}{\mathcal{M}}
\renewcommand{\S}{\mathcal{S}}

\newcommand{\T}{\mathcal{T}}
\newcommand{\A}{\mathfrak{A}}

\newcommand{\lb}{\left (}
\newcommand{\br}[1]{\left( #1 \right ) }

\newcommand{\rb}{\right )}
\newcommand{\su}[1]{\mathop{\sum}_{#1}}
\newcommand{\htM}{\widehat{\mathcal{M}}}

\newcommand{\defn}[1]{\begin{defin} \  #1 \end{defin}}
\newcommand{\prop}[1]{\begin{propo} \  #1 \end{propo}}

\begin{document}
    \title{
\begin{flushright}
 \mbox{\normalsize ITEP/TH-97/09}
\end{flushright}
\vskip 20pt
On the Homology of Certain Smooth Covers of Moduli Spaces
of Algebraic Curves}
    \author{P.Dunin-Barkowski\footnote{ITEP, Moscow, Russia and KdVI, UvA, Amsterdam, The Netherlands; barkovs@itep.ru}, A.Popolitov\footnote{ITEP, Moscow, Russia; popolit@itep.ru}, G.Shabat\footnote{ITEP and RSHU, Moscow, Russia; george.shabat@gmail.com}, A.Sleptsov\footnote{ITEP, Moscow, Russia; sleptsov@itep.ru}}
    \maketitle
\begin{abstract}
We suggest a general method of computation of the homology of certain smooth covers $\widehat{\mathcal{M}}_{g,1}(\mathbb{C})$ of moduli spaces $\mathcal{M}_{g,1}\br{\mathbb{C}}$ of pointed curves of genus $g$. Namely, we consider moduli spaces of algebraic curves with level $m$ structures. The method is based on the lifting of the Strebel-Penner stratification $\mathcal{M}_{g,1}\br{\mathbb{C}}$. We apply this method for $g\leq 2$ and obtain Betti numbers; these results are consistent with Penner and Harer-Zagier results on Euler characteristics.
\end{abstract}

\tableofcontents

    \section {Introduction}
The homology of moduli spaces of curves deserves much attention
during the last decades ( see, e.g., \cite{F1, F2}). However, the values
of \it all \rm the Betti numbers $b_q(\mathcal{M}_g)$ are far from
being known. There exist some  indications of the existence of the
\it beautiful \rm answer for these numbers; the values of orbifold
Euler characteristics \cite{HZ} and \cite{P} and the generating functions for the intersection numbers \cite{K} and \cite{W} are among the best known.

Most results on the homology of moduli spaces were obtained by
some indirect methods; the papers by Looijenga \cite{L1,L2} is among the few
counterexamples known to the authors. 

In the present paper we take the direct approach to the calculation
of $\mathrm{H}_*(\mathcal{M}_g,\mathbb{Z}/2\mathbb{Z})$ based on the stratification of
the moduli spaces labeled by \it dessins d'enfants. \rm 
The two
versions of these stratifications were introduced by Penner in \cite{P2} (see also \cite{P3})
and by Strebel (see Looijenga \cite{L1}); they are different set-theoretically but equivalent
combinatorially.

The similar approach was undertaken in \cite{PMil}. 

In order to avoid difficulties related to singularities of moduli space we work with smooth cover instead and lift the Strebel-Penner decomposition there. 
Then we construct a simplicial complex on which this cover retracts, and to which the standard definition of simplicial homology can be applied.

Another motivation of our work comes from the fact, that moduli spaces $\mathcal{M}_{g,\gamma}$ play important role in string theory. Namely, in order to calculate string scattering amplitude of $\gamma$ (closed) strings, one should perform an integration over $\mathcal{M}_{g,\gamma}$ (this is analogous to the integration over momenta of virtual particles in Feynmann diagram technique for particle physics). 
There are many open problems in this field, even the measure of integration in supersymmetric case is known only up to genus 4 (for developement in this field see \cite{DP,DBMS}).

Our paper is organized as follows. In section 2 we show that using dessins d'enfants one defines a structure of cellular complex on $\mathcal{M}_{g,\gamma}\times\mathbb{R}_{>0}^\gamma$. For this reason we first define metrized ribbon graphs as pairs of dessins d'enfant and real positive numbers. Then we homeomorphically map metrised ribbon graphs to the moduli space that provides us a decomposition of moduli space into cells. Then in section 3 we construct a smooth cover $\widehat{\mathcal{M}}_{g,1}$ of $\mathcal{M}_{g,1}$ in order to resolve the singularities of moduli spaces. It is done by introducing an additional structure: a basis in homology of Riemann surfaces. In section 4 we construct a spine $\mathcal{S}$ on $\widehat{\mathcal{M}}_{g,1}$ and prove that it is a simplicial complex. Then we conjecture that the simplicial complex can be represented as a cellular complex $\mathbb{D}$ by combining groups of simplicies into cells. In section 5 we describe our method for calculating the Betti numbers of $\widehat{\mathcal{M}}_{g,1}$. For simplicity it is divided into practical steps. In sections 6 and 7 our method is applied to genus 1 and 2. Explicit answers are given.

\section{Metrized ribbon graphs and moduli spaces}
\subsection{ Dessins d'enfants} 
We briefly introduce the main concepts and terminology of the
theory (see \cite{VS} and \cite{LZ} for more details). \\
\\
A triple
$$
X_0\subset X_1\subset X_2=X
$$
is called a \it dessin d'enfant \rm if\\
\\
$\bullet X_0$ is a finite set (of \it vertices \rm);\\
$\bullet X_1$ is a graph, i.e.
$$
X_1\setminus X_0\text{ is homeomorphic to a disjoint union of finite
number of open real intervals (\it{edges})};
$$
 $\bullet X_2$ is a compact connected oriented surface and
 $$
X_2\setminus X_1\text{ is homeomorphic to a disjoint union of finite
number of open discs (\it{cells})}.
$$\\
Throughout the paper we denote by $g$ the genus of a surface $X$ and
by $\gamma$ the number of the cells; we are going mostly to consider
the case $\gamma=1$. Throughout the paper we assume
$(g,\gamma)\notin\{(0,1),(0,2)\}$.\\
\\
Following \cite{Gr}, introduce the \it oriented cartographic group \rm$\mathcal{C}_2^+$; as
an abstract group it is generated by three elements $\rho_{0,1,2}$ satisfying the relations
$$
\rho_2\rho_1\rho_0=1
$$
and
$$
\rho_1^2=0.
$$
(Of course, using the first relation we can delete one of the
generators $\rho_0$ or $\rho_2$ but it is convenient to use both).\\
\\
The group \rm$\mathcal{C}_2^+$ should be thought of as acting on the
\it directed edges \rm of any dessin d'enfant; $\rho_0$ rotates a
directed edge contrary-clockwise around its origin (sending it to
the \it next \rm one), $\rho_1$ changes the direction and $\rho_2$
moves an edge contrary-clockwise around the cell that lies \it to
the left \rm of the directed edge. \\
\\
So the set $\vec{E}(D)$ of the directed edges of a dessin $D$ is a
homogeneous finite $\mathcal{C}_2^+$-set; it is well-known that a
dessin can be restored by this set.\\
\\
By a 0-\it valency \rm of a directed edge $\epsilon\in\vec{E}(D)$ we mean its
$\rho_0$-order, i.e. the minimal positive integer $k$ such that
$\rho_0^k\cdot\epsilon=\epsilon$; obviously this number equals the number of germs of
edges incident to the origin of $\epsilon$.\\
\\
For any pair $(g,\gamma)$  our main tool is the \it finite \rm set
$\mathrm{Dess}_{g;*,*,\gamma[v_0\ge3]}$ of the isomorphic classes of
dessins of genus $g$ with $\gamma$ cells and \it all \rm the
0-valencies exceeding 3. It can be thought of as finite number of
isomorphic classes of objects in the category of $\mathcal{C}_2^+\mathfrak{sets}$.\\
 \\
 For $\vec{E}\in\mathrm{Dess}_{g;*,*,\gamma}$ denote
 $$
 X_{\vec{E}}
 $$
 the corresponding compact surface of genus $g$ with $\gamma$ cells.
 \\
 \\

\subsection{ Metrised ribbon graphs}

 Given $\vec{E}\in\mathrm{Dess}_{g;*,*,\gamma[v_0\ge3]}$, consider the
 real octant
 $$
 \mathrm{Metr}_{\vec{E}}:=\{\ell\in\mathbb{R}_{>0}^{\vec{E}}|
 \rho_1^*\ell=\ell\}.
 $$
 \\
  A pair $(\vec{E},\ell)$ with
  $\vec{E}\in\mathrm{Dess}_{g;*,*,\gamma[v_0\ge3]}$ and
$\ell\in\mathrm{Metr}_{\vec{E}}$
  will be called a \emph{metrised ribbon graph}.\\
  \\
  To each
  $(g,\gamma)\in(\mathbb{N}\times\dot{\mathbb{N}})\setminus\{(0,1),(0,2)\}$
  we associate the finite disjoint union of real octants
  $$
  \mathrm{MRG}_{g,\gamma}:=
\coprod_{\vec{E}\in\mathrm{Dess}_{g;*,*,\gamma[v_0\ge3]}}\mathrm{Metr}_{\vec{E}}
   \cong\{(\vec{E},\ell)|\vec{E}\in\mathrm{Dess}_{g;*,*,\gamma[v_0\ge3]},
    \ell\in\mathrm{Metr}_{\vec{E}}\}.
  $$
  \\
  The most important are the $\rho_2$-invariant \emph{circumferences}
  $$
  c(\epsilon):=\sum_{\epsilon'\in<\rho_2>\epsilon}\ell(\epsilon').
  $$

  \subsection{Strebel atlases}  
To a metrised ribbon graph
  $(\vec{E},\ell)\in \mathrm{MRG}_{(g,\gamma)}$ we associate a
  covering
  $$
  X_{\vec{E}}=:\cup_{\epsilon\in\vec{E}}(\mathcal{Z}_\epsilon\cup
  \mathcal{U}_\epsilon\cup\mathcal{W}_\epsilon)
  $$
together with the local coordinates
$$
z_\epsilon:\mathcal{Z}_\epsilon\hookrightarrow\mathbb{C},\quad
u_\epsilon:\mathcal{U}_\epsilon\hookrightarrow\mathbb{C},\quad
w_\epsilon:\mathcal{W}_\epsilon\hookrightarrow\mathbb{C},
$$
\begin{figure}[H]
\centering\leavevmode
\includegraphics[width=15 cm]{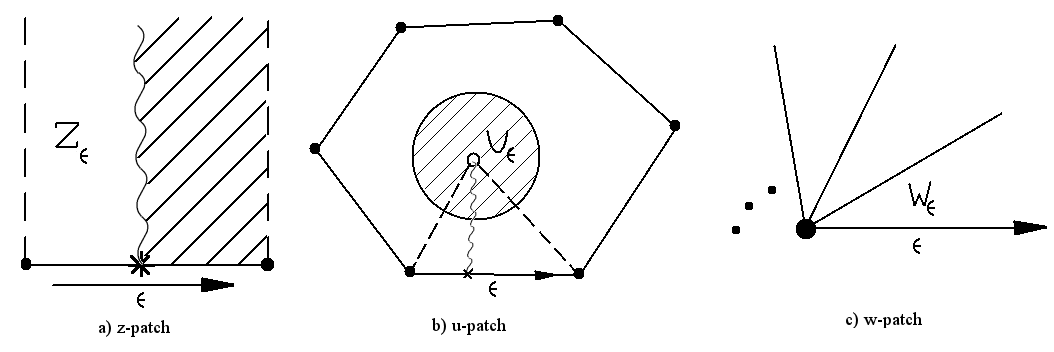}
\caption{Patches}
\label{basis}
\end{figure}
satisfying
 $$
 z_{\rho_1\cdot\epsilon}=-z_\epsilon,
 $$
 %%%%%%%%%%%
 $$
 z_{\rho_2\cdot\epsilon}=z_\epsilon+
 \frac{\ell(\epsilon)+\ell(\rho_2\cdot\epsilon)}{2};
$$
%%%%%%%%%%%%
 $$
 u_\epsilon=e^{\frac{2\pi\mathrm{i}}{c(\epsilon)}z_\epsilon};
 $$
%%%%%%%%%%%%
$$
w_\epsilon=[z_\epsilon+\frac{\ell(\epsilon)}{2}]^{\frac{2}{v_0(\epsilon)}}.
$$
\\
We have just defined a holomorphic \emph{atlas}
$$
\mathfrak{A}_{(\vec{E},\ell)}:=\{(\mathcal{Z}_\epsilon,z_\epsilon),
(\mathcal{U}_\epsilon,u_\epsilon),
(\mathcal{W}_\epsilon,w_\epsilon)|\epsilon\in\vec{E}\}
$$
on the topological surface $X_{\vec{E}}$.
\\

\subsection{From metrized ribbon graphs to moduli spaces}
The above construction defines a map
$$
\mathrm{MRG}_{g,\gamma}\longrightarrow\mathcal{M}_{g,\gamma}(\mathbb{C}).
$$
Compare the dimensions in the case of trivalent dessins $D$
($v_0\equiv3$):
$$
\dim_{\mathbb{R}}(\mathrm{MRG}_{g,\gamma})=\#\{\text{edges of }D\}=
 \frac{1}{2}\#\vec{E}(D);
$$
now, denoting $\alpha:=\#\{\text{vertices of }D\}$, taking into
account $\alpha-\frac{1}{2}\#\vec{E}(D)+\gamma=2-2g$ and
$\#\vec{E}(D)=3\alpha$ (that follows from our assumption
$v_0\equiv3$) and solving
$\alpha-\frac{1}{2}\#\vec{E}(D)+\gamma=2-2g$, where we substitute
$\alpha=\frac{1}{3}\#\vec{E}(D)$ we arrive at
$$
\dim_{\mathbb{R}}(\mathrm{MRG}_{g,\gamma})=3(2g-2+\gamma)=6g-6+3\gamma,
$$
while
$$
\dim_{\mathbb{R}}\mathcal{M}_{g,\gamma}(\mathbb{C})=
 2\dim_{\mathbb{R}}\mathcal{M}_{g,\gamma}(\mathbb{C})=2(3g-3+\gamma)=6g-6+2\gamma,
$$
so our map reduces real dimensions by $\gamma$. However, the points
of the  space $\mathrm{MRG}_{g,\gamma}$ define exactly $\gamma$
extra (positive) real parameters: the circumferences of the
boundaries of cells.\\
\\
The fantastic theorem of Strebel \cite{S} claims that
$$
\mathrm{MRG}_{g,\gamma}\longrightarrow\mathcal{M}_{g,\gamma}(\mathbb{C})
 \times\mathbb{R}_{>0}^\gamma,
$$
sending a metrized ribbon graph to the Riemann surface defined by
the Strebel atlas and a tuple of circumferences, actually
is a homeomorphism!

In what follows we set $\gamma=1$ and require the sum of the legths of all edges to be equal to $1$, thus dropping $\mathbb{R}_{>0}^\gamma$ factor.

Now we can define strata on the MRG side and map them to the moduli space via our homeomorphism.
Namely, a stratum is an equivalence class, where two MRG's belong to the same class iff they can be transformed
to one another by continiously changing lengths of the edges without letting any edge's length to go zero.
\begin{equation}
 \label{mstrat}
\mathcal{M}_{g,1}=:\mathop{\bigsqcup}_{D\in Dess}C_D
\end{equation}

\section{Covers $\widehat{\mathcal{M}}_{g,1}\longrightarrow
\mathcal{M}_{g,1}$}
\label{cover}
In this section we define $\widehat{\mathcal{M}}_{g,1}$, the so-called \textit{level $m$ covers of $\mathcal{M}_{g,1}$} \cite{ACG}.

\subsection{Definition}
Fix a finite abelian group $\mathbb{A}$.
Then by definition,
$$
\widehat{\mathcal{M}}_{g,\gamma}:=\{(\mathcal{X},\iota)\}
$$
with a fixed isomorphism
$$
\iota:\mathrm{H}^1(\mathcal{X},\mathbb{A})\mathop{\longrightarrow}^{\sim}
\mathbb { A } ^ { 2g } ;
$$
the map $\widehat{\mathcal{M}}_{g,\gamma}\longrightarrow
\mathcal{M}_{g,\gamma}$ forgets $\iota$.

In other words, $\widehat{\mathcal{M}}_{g,1}$ is a quotient of Teichmuller space $\mathcal{T}_{g,1}$ by a certain subgroup $\widehat{\Gamma}_{g,1}$ of finite index of modular group $\Gamma_{g,1}$, defined by the short exact sequence
\eq{
1 \longrightarrow \widehat{\Gamma}_{g,1}  \longrightarrow  \Gamma_{g,1}  \longrightarrow \mathrm{Sp}_{2g}\br{\mathbb{A}} \longrightarrow 1
}

\subsection{Smoothness}

For $\mathbb{A} = \mathbb{Z}/3\mathbb{Z}$ (as well as for all $\mathbb{A} = \mathbb{Z}/m\mathbb{Z}$, $m\geq 3$) the above defined $\htM_{g,1}$ are smooth (see Chapter XVI of \cite{ACG}). For $g=1$ the cover of level 2 is also smooth, so we consider this cover when we describe $g=1$ case in detail in Section \ref{gen1}.

\subsection{Stratification}

The stratification (\ref{mstrat}) can be naturally lifted to $\htM_{g,1}$. The strata are then enumerated by dessins d'enfants with a chosen symplectic basis in the first homology of the curve as an additional
structure, $\widehat{D} = (D,\iota)$.

In a particular such stratum, corresponding to a dessin with $k$ edges, all
points are uniquely described by $k$
numbers, sum of which is equal to $1$, i.e. the highest dimension of a stratum
is $6g-4$, because the maximal number of edges in our dessins is $6g-3$ (while
the minimal is $2g$).

It is evident that every stratum is in fact an (open) simplex. The stratification,
however, is \textbf{not} a simplicial complex, because not all faces of each
simplex lie in this stratification. This happens because for some dessins one
cannot retract some of their edges without pinching a handle or a throat. Let us call this stratification a $quasisimplicial complex$, meaning a difference of two simplicial complexes. Denote it as $\mathcal{Q}_g$, i.e.
\eq{
\htM_{g,1} \cong \mathcal{Q}_g = \mathop{\bigsqcup}_{\widehat{D}} C_{\widehat{D}}
}

\section{Spine}

\subsection{Definition}
\label{spine}

Consider a \textit{flag} $F$ of dessins d'enfants with additional structures,
\eq{
F : \widehat{D}_0 \prec \dots \prec \widehat{D}_k,
}
where $\widehat{D} \prec \widehat{D}'$ means that $\widehat{D}'$ can be obtained from $\widehat{D}$ by contraction of some of its edges.

Denote the set of all flags by $\mathcal{F}$.

Recall that to each dessin with additional structure $\widehat{D}$ corresponds an open simplex $S_{\widehat{D}}$ in quasisimplical complex $\mathcal{Q}_{g}$. Note that for all $i$ the stratum $S_{\widehat{D}_i}$ lies in the closure of $S_{\widehat{D}_k}$.

Denote by $S_F$ the open convex hull of all centers of $S_{\widehat{D}_i}$ in the coordinates of the closure of $S_{\widehat{D}_k}$.

Here by ``the center of a stratum'' we mean the point which corresponds to all edges of the dessin being of equal length. Note that $S_F$ is a $k$-simplex due to all of its vertices being in general position by construction.
\defn{
Define the spine of $\hM$ as
\eq{
\S_{g,1} := \bigsqcup_{F\in \Fl} S_F
}
}
\prop{
$\S_{g,1} $ is a simplicial complex.
}

By construction, $\mathcal{F}$ is an abstract simplicial complex with dessins with additional structures as its vertices and flags as faces. Since simplices of our quasisimplicial complex $\mathcal{Q}_g$ are glued together according to this very partial ordering on dessins with additional structures on which flags are based, the geometric realization of $\mathcal{F}$ with simplices $S_F$ is a simplical complex.

\subsection{Retraction}

\prop{

$\S_{g,1}$
is a deformation retract of $\widehat{\mathcal{M}}_{g,1}$.
}

To prove this proposition we will construct a continuous function $f: \hM
\times [0,1] \rightarrow \hM$, such that $f_0$ is identical on $\hM$ and
$f_1(\hM)=\S_{g,1}$.

First, we will define $f$ on a certain dense subset of $\hM \times [0,1]$ and then
we will prove that it can be extended to all $\hM \times [0,1]$ by continuity.

Let $S_{\widehat{D}_0}$ be a stratum of the highest dimension (recall that it is actually
an open simplex, as all strata are). Introduce a particular ordering on the edges $\widehat{D}_0$. 
Then we have coordinates $(l_1,\dots,l_{N}),\; \su{i}
l_i = 1$, on the stratum. Consider a subset of the stratum with
$l_1>l_2>\dots>l_{N}$. Denote it by $\A$. We will now define the
function $f$ on $\A\times [0,1]$. Note that by choosing all other possible
orderings of edges and considering analogous subsets, one covers all stratum.
Thus, if we define $f$ on $\A$ and then define it in analogious way on other
such subsets corresponding to other orderings and on other strata of the highest
dimension, we will have $f$ defined on a dense subset of $\hM \times [0,1]$.

Let
\eq{L = (l_1,\dots,l_{N})}
be a point in $\A$. Let us retract the edge of the
dessin which corresponds to the last coordinate (which is the smallest one by
construction of $\A$). Denote the resulting dessin as $\widehat{D}_1$. Then let us try to retract the edge which corresponds to
the next to last coordinate and so on. We will finally face a situation when we
try to retract an edge, retraction of which inevitably leads to pinching a
handle. Let $k+1$ be the index of the coordinate for which we faced such a
situation for the first time. Then denote by $O$ the following point:
\eq{
O = \dfrac{1}{l_1+\dots+l_{k} -
k l_{k+1}}\br{l_1-l_{k+1},l_2-l_{k+1},\dots,l_{k}-l_{k+1},0,\dots,0}
}
This point does not belong to $S_{\widehat{D}}$. Actually, it does not belong to $\hM$.
It lies on a boundary of $S_{\widehat{D}}$ not included in $\hM$.
Consider then the ray $OL$ in $S_{\widehat{D}}$. Note that it intersects one of
simplices of $\S_{g,1}$, namely the simplex corresponding to flag
$F: \widehat{D}_0 \prec \dots \prec \widehat{D}_k$.

To see it, recall the construction
of this simplex. Since this simplex is the convex hull of centers of $S_{\widehat{D}_0},\dots,S_{\widehat{D}_0}$, in this simplex lie all points of $S_{\widehat{D}_0}$ with first $k$ coordinates equal to each other and last $N-k$ ones strictly decreasing and less then the first ones.

It is straightforward to show that ray $OL$ passes through the point
\eq{
P(L) =
\dfrac{1}{1+k l_k-l_1-\dots-l_{k}}\br{l_{k+1},\,l_{k+1},\dots,l_{k+1},\,l_{k+2},\dots,
l_N},
}
where $l_{k+1}$ stand on first $k=1$ places and $N=6g-3$ is the number of edges of $\widehat{D}_0$. Obviously, this point lies on $OL$ further from $O$ then $L$ does.

Because $l_k>l_{k-1}>\dots>l_{N}$ it is evident that this point lies in $S_F$ due to all mentioned requirements being satisfied. It is also evident that this is the only point of intersection of $OL$ with $S_F$, because this is the only point on $OL$, where first $k$ coordinates are equal to each other (there cannot be more such points because $OL$ contains a point, namely point $L$, which has all coordinates mutually \textit{inequal}).

Let then $f_t$ map the point $L$ to the point on $LP$ which divides this segment as $t:\br{1-t}$, i.e.
\eq{
f:\br{L,t} \mapsto L-\br{L-P(L)}t
}
Because this map is linear in $t$ and all $l_i$, it is evidently continuous on $\A$. Also it is evident that $f_0$ is identical and that $f_1$ maps $\A$ on $S_F$ (one only needs to check that something is mapped to \textit{every} point of $S_F$; this can be done by just considering a point in $S_F$ and a point on the corresponding boundary and than taking any point on the connecting interval).

Thus we have defined a continuous map on a dense subset of $\hM\times [0,1]$. 
It is easy to see that this map can be continued to entire $\hM$.

\subsection{``Dual'' cell complex}
\label{cells}
All information about homology of $\widehat{\mathcal{M}}_{g,1}$ is encoded in its spine $\S_{g,1}$, a pure simplicial complex of dimension $4g-3$, defined and studied above. However, already for genus 2 the number of simplices in this simplicial complex is far too huge for its homologies being computable on current hardware. Thus we prove that actually it can be represented as a \textit{cellular complex} by combining collections of simplices into cells, with eventual number of cells being far smaller than initial number of simplices.

This cell complex is constructed in the following way.

Let $\T_0$ be the set of centers of all strata of codimension $4g-3$. Evidently, $\T_0 \subset \S_{(0)}$, i.e. $\T_{0}$ is a certain subset of the set of vertices of the simplicial complex. For $T \in \T_{0}$ consider its star $A_T$ in $\S_{g,1}$ (recall that the star of a vertex of a simplicial complex is the union of all simplices of this complex which have this vertex as one of their vertices). Denote the set of all such stars as $\As_{(4g-3)}$. Exclude all simplices that formed these stars from $\S$. The remaining set of simplices will still be a simplicial complex.

Consider then the set $\T_1$, centers of all strata of codimension $4g-4$. Consider stars of elements of this set in the mentioned remaining simplicial complex. They form the set $\As_{(4g-4)}$ and so on.

\begin{lem}
\label{celllemma}
$\As = \mathop{\cup}_{i=0}^{4g-3}\As_{(i)}$ is a cell complex.
\end{lem}
In order to prove this, let us first prove the following proposition.

\prop{
Stars $A_{(4g-3)}(T),\ T \in \T_{0}$ as simplicial complexes are isomorphic to barycentric subdivision of \textit{associahedron $\mathcal{K}^{(4g-3)}$}.
}
For the definition of $\mathcal{K}^{(4g-3)}$ refer to e.g. \cite{Loday}.

Let $T$ be an element of $\T_{0}$ and let $\hD_0$ be the dessin with structure corresponding to the stratum of $\hM$ where $T$ belongs. Note that simplices of its star $A_{(4g-3)}(T)$ correspond to all possible flags of dessins with structures such that $\hD_0$ can be obtained from them by contraction of some number of edges. Let $\mathcal{D}_T$ be the set of all such dessins with structures and let $\mathcal{F}_T$ be the set of all such flags ($A_T$ is then a realization of $\mathcal{F}_T$ as an abstract simplicial complex).

Note the following 
\prop{ \label{propdes} For a given $\hD$ and two sets $E_1$ and $E_2$ of its edges, if $E_1$ differs from $E_2$, then dessins with structures $\hD_1$ and $\hD_2$ obtained from $\hD$ by contraction of all edges in $E_1$ and  $E_2$ respectively are different.} 
This proposition follows from the fact that a curve with an additional structure of this type does not have any nontrivial automorphisms \cite{ACG}. Let us stress that this proposition is not valid for dessins without additional structures, but it \textit{is} true for dessins with chosen symplectic bases.

This fact leads to the following statement: $\mathcal{D}_T$ is in one-to-one correspondance with the set of all planar trees with $4g-1$ leaves. If one cuts all of the edges of the original dessin $\hD_0$, one will obtain a planar tree. All other dessins in $\mathcal{D}_T$ again turn into trees if one cuts these $4g$ edges in them. Then the previous fact implies that it is a one-to-one correspondance.

These trees label the cells of associahedron $\mathcal{K}^{4g-3}$ (which is a cell complex). A tree with $k$ internal vertices corresponds to a $(4g-2-k)$-cell. The barycentric subdivision of $\mathcal{K}^{4g-3}$ is a realization of an abstract simplicial complex formed by flags of trees. These flags are actually the same as the flags in $\mathcal{F}_T$, via the mentioned correspondance. This implies that $A_T$ is isomorphic to the barycentric subdivision of $\mathcal{K}^{4g-3}$.

This proposition implies that all stars $A_T$ for $T \in T_0$ are homeomorphic to balls of corresponding dimension, due to the associahedron being homeomorphic to a ball.

Analogously, all defined above stars of lower dimesion are homeomorphic to balls as well, which implies the statement of the lemma.

Note that resulting cell complex is ``dual'' to quasisimplicial complex $\mathcal{Q}_g$ in the sense that simplices of the highest dimension of $\mathcal{Q}_g$ are in one-to-one correspondance with vertices of cell complex, simplices of codimension one are in one-to-one correspondance with edges of cell complex and so on.

Also note that the incidence number arising in computation of cellular homology are trivial in the case of complex $\mathcal{A}$. This again follows from Proposition \ref{propdes}.

\section {Description of the computational method}
Let us briefly repeat the practical steps needed to be taken in order to actually compute Betti numbers of the moduli spaces of level $m$ of algebraic curves of genus $g$, where $m=2$ for $g=1$ and $m=3$ otherwise. In the two sections that follow this one we present the details of computation for $g=1$ and $g=2$.

In this practical approach we actually compute the homology of the cell complex $\mathcal{A}_g$ which is a deformation retract of $\hM$. We use the fact that this complex is ``dual'' to $\mathcal{Q}_g$.

The practical steps are the following ones:
        \begin{enumerate}
            \item Describe all possible dessins d'enfants with at least trivalent vertices that can be drawn on a surface
            of a given genus , which give exactly one disc when cut along. It is equivalent to enumerating
            all inequivalent ways of glueing of a $(12g-6)$-gon.
            \item For each glueing scheme find all symplectic bases in 1-homology of the curve. For
            this one needs to define an intersection form on the glueing scheme (see Section \ref{intersect} for details).
            \item Construct the structure of the cell complex by labeling cells with gluing scheme with basis.
            \item Define the boundary operator corresponding to contraction of edges step by step.
            \item Compute dimensions of kernels and images of the boundary operator.
        \end{enumerate}

        \section{Genus 1 case}
\label{gen1}
In this section we compute Betti numbers of $\htM_{1,1}$, where $\htM_{1,1}$ is a six-sheeted covering of $\Mt_{1,1}$.

First, let us enumerate all needed graphs.  For this we calculate Euler's characteristic and impose condition, that graph is trivalent:
\eq{
\chi \equiv V - E + F = 2 - 2g = 0
}
\eq{\left\{
\begin{aligned}
& V - E + 1 = 0 \hspace{0.4cm} \Leftrightarrow \hspace{0.3cm} \\
& 3V = 2E \\
\end{aligned}
\right.
\left\{
\begin{aligned}
V =2 \\
E=3 \\
\end{aligned}
\right.}

There is a unique graph, that satisfies these requirements, which is called $\theta$.

\begin{figure}[!ht]
\begin{center}
\includegraphics[scale=0.3]{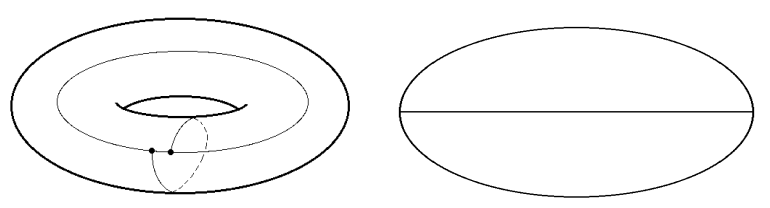}
\end{center}
\caption{$\theta$-graph}
\label{pic}
\end{figure}

Second, we define functions on the edges of the graph over $\mathbb{Z}/2\mathbb{Z}$, i.e. we associate $0$ or $1$ with each edge, so that the sum at each vertex is equal to zero. Then there are $4$ functions on this graph:
\eq {
\lb \begin{array}{c}
0  \\
0  \\
0
\end{array} \rb, \
e_1:=\lb \begin{array}{c}
1  \\
1 \\
0
\end{array} \rb, \
e_2:=\lb \begin{array}{c}
0  \\
1 \\
1
\end{array} \rb, \
e_3:=\lb \begin{array}{c}
1  \\
0  \\
1
\end{array} \rb.
}
We consider these functions as vectors generating $2$-dimensional vector space. To come to agreement with canonical intersection form of cycles the basis of the vector space has to be symplectic. Let us enumerate all symplectic bases:
\eq{
(e_1,e_2), \ (e_1,e_3), \ (e_2,e_3), \ (e_2,e_1), \ (e_3,e_1), \ (e_3,e_2).
}

Third, glueing schemes with symlectic bases form cell. Since in $g=1$ case there is the only glueing scheme, then there are $6$ cells as well as a number of symplectic bases.

Let us introduce the following notation: if the edge is labeled by $\br{0,1}$, we denote this edge by $\textbf{a}$;  if the edge is labeled by $\br{1,0}$, we denote it by $\textbf{b}$, otherwise $\textbf{c}$. All zeroes situation is not considered since it is not presented in any of abovewritten bases superficially. So, there are 6 cells of highest dimension in our complex:

\begin{figure}[H]
\centering\leavevmode
\includegraphics[width=12 cm]{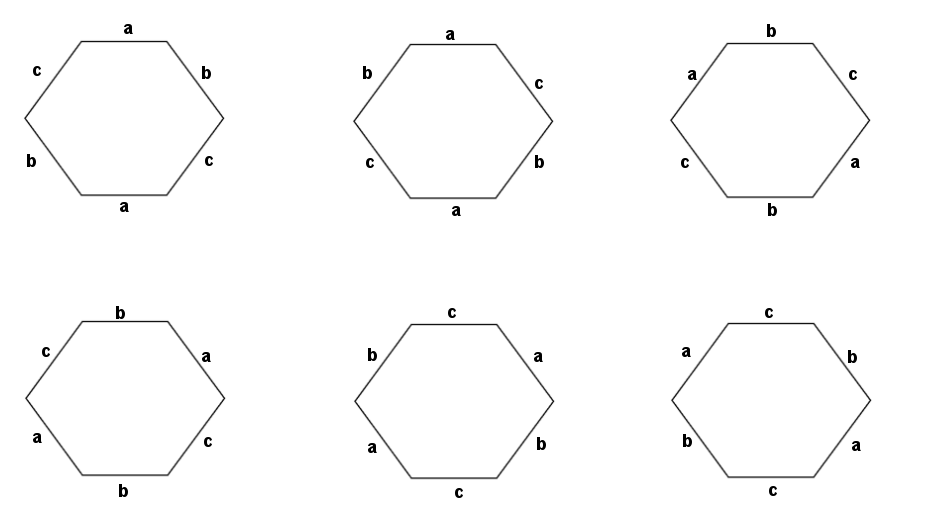}
\caption{FULL BASIS}
\label{Fbasis}
\end{figure}

Note that some glueing schemes with bases are equivalent with respect to rotation. Thus, there are only 2 nonequivalent glueing schemes with bases, i.e. we have $2$ cells of highest dimensions:
\begin{figure}[H]
\centering\leavevmode
\includegraphics[width=9 cm]{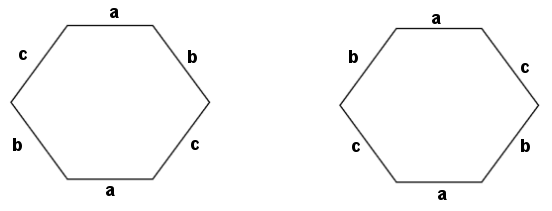}
\caption{REDUCED BASIS}
\label{rbasis}
\end{figure}

Fourth, to find out how these cells are connected to each other, we contract each edge of $\theta$-graph by turn and again factorize by rotation. Thus, we have 3 nonequivalent glueing schemes:
\begin{figure}[H]
\centering\leavevmode
\includegraphics[width=12 cm]{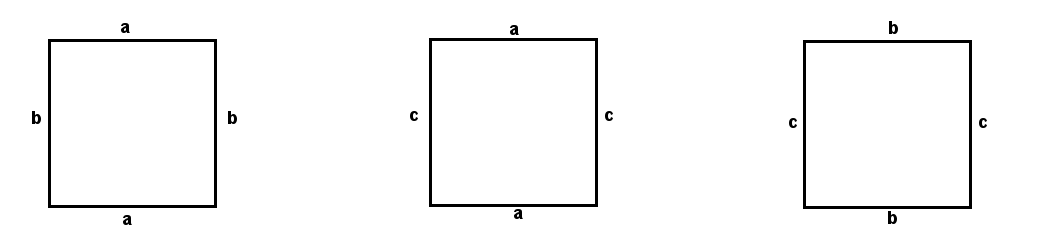}
\caption{1-CELLS}
\label{1cells}
\end{figure}

If we contract one more edge, than we will contract torus handle and genus will get down. Thus, three 0-cells $\textbf{a},\textbf{b},\textbf{c}$ are not included in our cellular complex. So, our cellular complex consists of two $2$-dimensional cells $\textbf{abc}$ and $\textbf{acb}$, three $1$-dimensional cells $\textbf{ab}$, $\textbf{bc}$ and $\textbf{ac}$. Geometrically, we can illustrate the complex as follows:
\begin{figure}[H]
\centering\leavevmode
\includegraphics[width=5 cm]{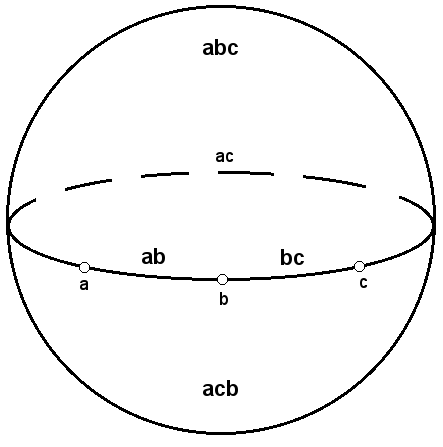}
\caption{CELLULAR COMPLEX}
\label{complex}
\end{figure}

Fifth, we construct dual complex by the following way $$ {\rm 2-}cell \ \longmapsto \ {\rm 0-}cell \hspace{2cm} {\rm1-}cell \ \longmapsto \ {\rm 1-}cell $$
\begin{figure}[H]
\centering\leavevmode
\includegraphics[width=8 cm]{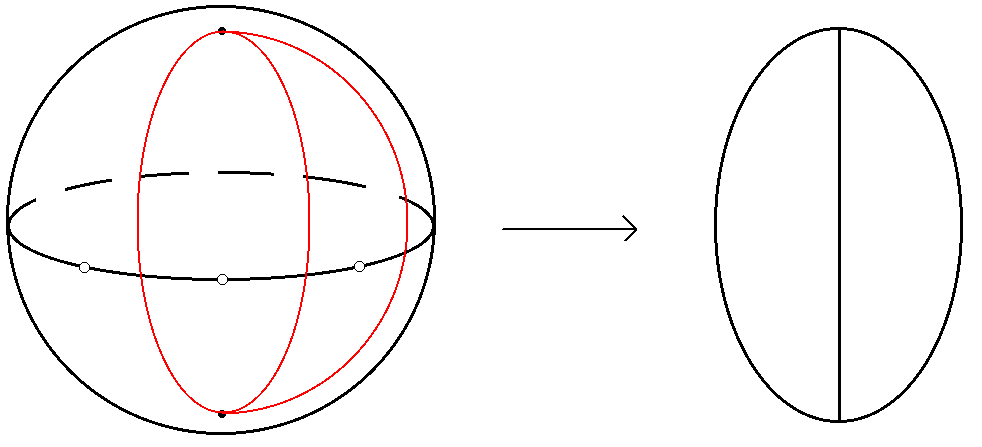}
\caption{DUAL CELLULAR COMPLEX}
\label{ducomplex}
\end{figure}

Thus, clearly, Betty numbers of moduli space for genus 1 with 1 marked point are
\eq{
\begin{array}{|c|}
\hline\\
dim\mathbb{H}^0(\Mt_{1,1},\mathbb{Z}/2\mathbb{Z})=1\\
dim\mathbb{H}^1(\Mt_{1,1},\mathbb{Z}/2\mathbb{Z})=2\\
dim\mathbb{H}^k(\Mt_{1,1},\mathbb{Z}/2\mathbb{Z})=0, k \geq 2\\ \\
\hline
\end{array}
}

\section {Genus 2 case}

First of all, for this $g=2$ case we computed the homology of stars and it turned out to be trivial, thus allowing us to compute the homology of the moduli space as cellular homology of the quotient of spine by star structure rather than just as simplicial homology of the spine itself.

\par 
Also, the case of genus two is rather special,
due to existence of hyperelliptic involution the covering of $\Mm_{2,1}$
obtained with the help of symplectic bases in homologies of the surface over field
$\mathbb{Z}/2\mathbb{Z}$ is still not smooth and one needs to consider symplectic bases in homologies
over $\mathbb{Z}/3\mathbb{Z}$.

The results on ribbon graphs on genus-2 surfaces are summarised in the following table:

\begin{tabular}{|c||c@{\hspace{7mm}}|c@{\hspace{7mm}}|c@{\hspace{7mm}}|c@{\hspace{7mm}}|c@{\hspace{7mm}}|c@{\hspace{7mm}}|c@{\hspace{7mm}}|c@{\hspace{7mm}}|c|}
 \hline
 number &  \multicolumn{8}{|c|}{number of schemes with symmetry of order  } &  total number \\
 
\cline{2-9}
 of edges & 1  & 2 & 3 & 4 & 5 & 6 & 8 & 10 & of schemes \\
\hline
\hline
9 & 3 & 5 & 1 &   &   &   &   &    & 9 \\
\hline
8 & 24 & 4 &   & 1 &   &   &   &    & 29 \\
\hline
7 & 41 & 11  &   &  &   &   &   &    & 52 \\
\hline
6 & 37 & 5 & 1  & 1 &   & 1  &   &    & 45 \\
\hline
5 & 14 & 5 &   & 4 & 1  &   &   & 1  & 21 \\
\hline
4 & 2 & 1 &   &  &   &   & 1  &    & 4 \\
\hline
\end{tabular}

\

On an asymmetric graph there are 51840 different symplectic bases in homologies over $\mathbb{Z}/3\mathbb{Z}$ on the curve. On a graph with a symmetry of order $r$ there are exactly $51840/r$ different bases. Graphs with chosen bases enumerate cells in our complex, with graphs with 9 edges corresponding to 0-cells in (dual) complex, with 8 edges -- to 1-cells, and so on. Thus, for the numbers $n_j$ of cells of dimension $j$ we have
\eqna{
n_0 = 51840\cdot & \displaystyle\br{\frac{1}{3}+\frac{5}{2}+3} &= 302400\\
n_1 = 51840\cdot & \displaystyle\br{\frac{1}{4}+\frac{4}{2}+24} &= 1360800\\
n_2 = 51840\cdot & \displaystyle\br{\frac{11}{2}+41} &= 2410560\\
n_3 = 51840\cdot & \displaystyle\br{\frac{1}{6}+\frac{1}{4}+\frac{1}{3}+\frac{5}{2}+37} &= 2086560\\
n_4 = 51840\cdot & \displaystyle\br{\frac{1}{10}+\frac{1}{5}+\frac{5}{2}+14} &= 870912\\
n_5 = 51840\cdot & \displaystyle\br{\frac{1}{8}+\frac{1}{2}+2} &= 136080\\
}
Euler characteristic of this complex is thus
\eq{
\chi_{\widehat{\Mm}_{2,1}} = 302400 - 1360800 + 2410560 - 2086560 + 870912 - 136080 = 432
}
Since Euler characteristic of a covering is equal to Euler characteristic of the base times the degree, we have
\eq{
\chi_{{\Mm}_{2,1}} = \displaystyle\frac{432}{51840} = \frac{1}{120}= \zeta\br{-3}
}
i.e. Euler characteristic of original noncompactified moduli space ${\Mm}_{2,1}$ is precisely equal to the value of Riemann zeta function at $-3$, which is in perfect agreement with results of Harer and Zagier.

As mentioned above, to compute Betti numbers of this covering of the moduli space it is suffecient to consider homologies of this cell complex over $\mathbb{Z}/2\mathbb{Z}$. We thus have 5 matrices filled with zeros and ones as matrices of boundary operators. Their ranks are
\eqna{
r_{1\rightarrow 0} &=& 302399\\
r_{2\rightarrow 1} &=& 1058377\\ 
r_{3\rightarrow 2} &=& 1352011\\
r_{4\rightarrow 3} &=& 733978\\
r_{5\rightarrow 4} &=& 135716\\
}
Thus, for Betti numbers we have
\eqna{
b_0 &=& 1\\
b_1 &=& 24\\ 
b_2 &=& 172\\
b_3 &=& 571\\
b_4 &=& 1218\\
b_5 &=& 364\\
}

\section{Conclusion}

In the present paper we constructed a direct method of computation of homology of certain smooth covers of moduli spaces of pointed curves via construction of the spine, a simplicial complex, on which the cover of moduli space retracts. Then we conjectured the equivalence of this spine to certain cell complex, which allows direct computation for the case of genus 2. We carried out this calculation and obtained Betti numbers for this case.

The proof of the mentioned conjecture remains a work in progress.

Since in our work we computed Betti numbers, then, recalling the mentioned in introduction connection with string theory, it seems natural to ask if there exist some generalized Penner model, i.e. a matrix model (or, perhaps, $\beta$-ensemble), that generates Betti numbers and not just Euler characterstics. It is very interesting if this model is an integral representation of some important quantity (such as Nekrasov function) on the string theory side. We are currently working in these directions.

\appendix
    \section {Details of computation}
        \subsection {Intersection form}
\label{intersect}
        
        \par
        In order to find possible symplectic bases on a given curve (parametrized by a glueing scheme),
        one needs to know an intersection form on cycles, that is
        bilinear map, which for any two cycles $a$ and $b$ gives their intersection number.
        Intersection form for $g = 1$ is very simple and is described in the corresponding section,
        so here we are focusing on the form for genus two. 
	Note, that while we consider curve homology with coefficients in $\mathbb{Z}/3\mathbb{Z}$ for $g=1$, for $g\geq 2$ we use homology with coefficients in $\mathbb{Z}/3\mathbb{Z}$.
        
        Note, that once we have found sympletic bases for glueing schemes of highest dimension, possible
        bases on cells of lower dimension are induced by contraction of some edges of cells of highest dimension.
        So, in what follows we concentrate on obtaining the intersection form on cells of highest dimension.

        \par Each one cycle, drawn on a curve, can be deformed in such a way, that it will run
        along edges of the glueing scheme, so each cycle is a map from edges of a glueing scheme to $\mathbb{Z}/3\mathbb{Z}$.
        Now, the answer is as follows: for any two cycles $a$ and $b$ their intersection number is equal to
        \eq {
            Int_{ab} = \lb \sum_{i} \Omega (a(i) \times a(i-1), b(i) \times b(i-1)) \rb \mod 3,
		\label {int_number}
	  }
        where sum is taken over vertices of a glueing scheme, $a(i)$ denotes the value
        of $a$ on the i-th edge. Note, that $\Omega$ defined below is neither \textit{antisymmetric}, nor
        \textit{polylinear}. Nevertheless, the intersection number, computed with help of it,
        is antisymmetric and polylinear function of $a$ and $b$.
        \par $\Omega$ is equal to
        \eq {
            \Omega = \lb \begin{array}{cccccccccc}
                          & 0\times0 & 0\times1 & 0\times-1 & 1\times0 & 1\times1 &
                          1\times-1 & -1\times0 & -1\times1 & -1\times-1 \\
                          0\times0: & 0 & 0 & 0 & 0 & 0 & 0 & 0 & 0 & 0 \\
                          0\times1: & 0 & 0 & 0 & 0 & 0 & 0 & 0 & 0 & 0 \\
                          0\times-1: & 0 & 0 & 0 & 0 & 0 & 0 & 0 & 0 & 0 \\
                          1\times0: & 0 & 0 & 0 & 0 & 0 & 0 & 0 & 0 & 0 \\
                          1\times1: & 0 & -\frac{1}{2} & \frac{1}{2} & \frac{1}{2} & 0 & 1 & -\frac{1}{2} & -1 & 0 \\
                          1\times-1: & 0 & \frac{1}{3} & -\frac{1}{3} & \frac{1}{3} & -\frac{1}{3} & -\frac{1}{2} & -\frac{1}{3} & \frac{1}{2} & \frac{1}{3} \\
                          -1\times0: & 0 & 0 & 0 & 0 & 0 & 0 & 0 & 0 & 0 \\
                          -1\times1: & 0 & -\frac{1}{3} & \frac{1}{3} & -\frac{1}{3} & \frac{1}{3} & \frac{1}{2} & \frac{1}{3} & -\frac{1}{2} & -\frac{1}{3} \\
                          -1\times-1: & 0 & \frac{1}{2} & -\frac{1}{2} & -\frac{1}{2} & 0 & -1 & \frac{1}{2} & 1 & 0                      
                     \end{array}
                    \rb
        }
        Let us turn to the derivation of this formula. Each cycle over $\mathbb{Z}/3\mathbb{Z}$ can
        be thought of as oriented ribbon, where $+1$ means that orientation of the ribbon coincides
        with that of an edge and $-1$ means that it is opposite. Also, since $1+1+1 = 0$, trivalent glueings
        of the ribbon are allowed, such that orientation of all glued ribbons are incoming (outgoing),
        respectively.
        \par Overall sign of the intersection number is a matter of convention.
        We use the following convention: left picture counts as $+1$ while right picture counts
        as $-1$.
        
\begin{figure}[!ht]
\begin{center}
\includegraphics[scale=0.5]{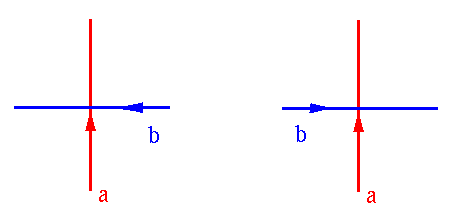}
\end{center}
\caption{Positive and negative intersections}
\label{pos_neg_int}
\end{figure}

        If one looks at some vertex for a given cycle, one can see one of the following pictures.
\begin{figure}[!ht]
\begin{center}
\includegraphics[scale=0.25]{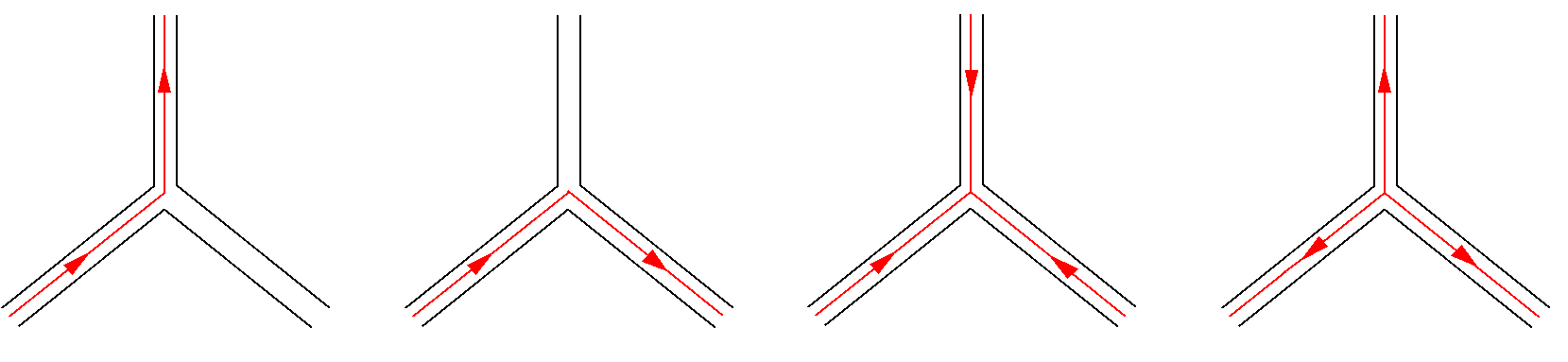}
\end{center}
\caption{Possible behavior of a cycles at the vertex}
\label{tri_vert_exampl}
\end{figure}        
        The idea is, that for each vertex, for each pair of possible pictures, an intersection
        number can be explicitly written.
        \par
        To succeed, without loss of generality assume, that whenever values of $a$ and $b$ on a given edge
        attached to some vertex are both non-zero $a$ is on the ``left'' side of the edge, while
        $b$ is on the ``right'' side:
\begin{figure}[!ht]
\begin{center}
\includegraphics[scale=0.4]{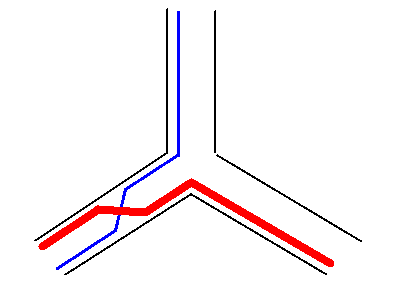}
\end{center}
\caption{ordering of cycles}
\label{cycle_order}
\end{figure}

\par One further distinguishes two kinds of situations: when both cycles are of type $L$ or $R$
at a given vertex and follow the same road (see pic), and all other situations.

\begin{figure}[!ht]
\begin{center}
\includegraphics[scale=0.4]{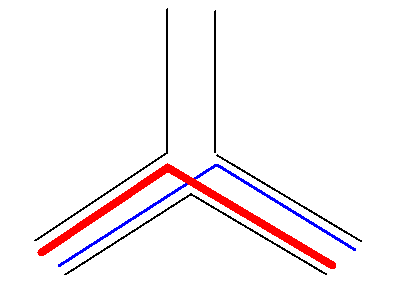}
\end{center}
\caption{Example of the trivial behavior at the vertex}
\label{trivial_vert_ex}
\end{figure}
Let us call situations (and vertices) of the first kind \textit {trivial}, others - \textit {nontrivial}.
It is easy to see that trivial vertices unify in sequences, in which $a$ and $b$ run along eash other,
and there are nontrivial vertices on both ends of such a sequence. Clearly, due to our ordering
convention, one should attribute either $+1$ or $-1$ intersection number to the whole sequence
of trivial vertices. Or, alternatively, one may prescribe an additional $\pm \frac{1}{2}$ to
nontrivial vertices, and do not consider trivial vertices at all.
\begin{figure}[!ht]
\begin{center}
\includegraphics[scale=0.4]{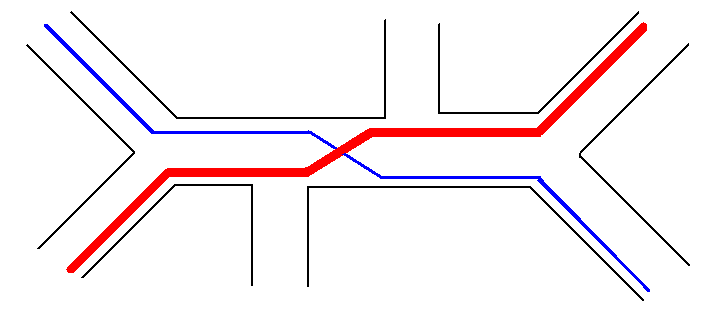}
\end{center}
\caption{Sequence of trivial vertices adds $+1$ or $-1$ to the intersection number}
\label{trivial_account}
\end{figure}

\par The next step is to consider each essential situation and prescribe intersection number
to it, which is fairly straightforward.
\par Next, one should notice, that type of the situation that actually takes place at given
vertex of a graph can be unambiguously reconstructed from values of cycles on just two of three
incoming edges (that is, by looking on vertex of a glueing scheme and edges, attached to it).
\par This way, one can rewrite sum over vertices of a graph on a curve as sum over vertices of
a glueing scheme. The only thing that one should take into account is, that risks to count
some vertices more than once, so symmetry factors like $\frac{1}{2}$ and ${1}{3}$ should be
introduced.
\par If one proceeds through all these steps, one arrives at the formula (\ref{int_number}).

\subsection {Calculation of ranks}
Matrices of the boundary operator, that appear in our calculation, are pretty large, so it is challenging to calculate their ranks. Luckily, there is a wonderful collection of linear algebra tools called\\
\textit{linbox} (http://www.linalg.org/linbox-html/index.html). So, we used the implementation of the Wiedemann algorithm, presented there, to calculate ranks of our matrices.

\section*{Acknowledgements}
Authors are very grateful to Sergey Shadrin for the idea of proof of Lemma \ref{celllemma} and to Elena Kreines for helpful discussions.
It is a special pleasure to thank Robert Penner for stimulating discussions and very helpful remarks. Authors are also grateful to Jim Milgram for a useful reference and to Jean-Guillaume Dumas for the help with \textit{linbox} software. This work is partly supported by NWO grant 613.001.021 (P.D.-B.); by RFBR grants 10-01-00709 (G.Sh.), 10-02-00509 (A.P., A.S.) and 10-02-00499 (P.D.-B.); by Federal Agency for Science and Innovations of Russian Federation (contract 02.740.11.0608); by Russian Federal Nuclear Energy Agency, by the joint grants  10-02-92109-Yaf-a (P.D.-B., A.P.), 09-02-91005-ANF (P.D.-B., A.P.), 09-02-93105-CNRSL (A.P.) and 09-01-92440-CE (P.D.-B., A.S.).

\end{document}